\documentclass{ifacconf}

\usepackage{graphicx}      
\usepackage{natbib}        

\usepackage{amsmath, amsfonts}
\usepackage{scalerel}
\usepackage{mathtools}

\newcommand{\mmin}{_{\min}}
\newcommand{\mmax}{_{\max}}
\newcommand{\pll}{^\mathrm{PLL}}
\newcommand{\cc}{^\mathrm{CC}}
\newcommand{\filter}{^\mathrm{f}}
\newcommand{\grid}{^\mathrm{g}}
\newcommand{\dqref}{_{dq}^\mathrm{r}}

\newcommand{\dlt}{\scaleobj{0.7}{\triangle}}

\newcommand{\dif}{\mathrm{d}}
\newcommand{\od}[2][]{\frac{\dif #1}{\dif #2}}

\newcommand{\const}{\mathrm{const}}

\newcommand{\bbR}{\mathbb{R}}
\newcommand{\bbZ}{\mathbb{Z}}

\newcommand{\setb}[1]{\bigl\{#1\bigr\}}

\newcommand{\p}[1]{(#1)}
\newcommand{\pb}[1]{\big( #1 \big)}

\newcommand{\norm}[1]{\Vert #1 \Vert}

\newtheorem{assumption}[thm]{Assumption}
\newtheorem{definition}[thm]{Definition}

\begin{document}
\begin{frontmatter}

\title{Numerical Estimation of the Lock-In Domain of a DC/AC Inverter\thanksref{footnoteinfo}}

\thanks[footnoteinfo]{The authors gratefully acknowledge funding by the German Federal Ministry of Education and Research (BMBF) within the Kopernikus Project ENSURE “New ENergy grid StructURes for the German Energiewende” (03SFK1B0-3).}

\author[First]{Anton Ponomarev}
\author[First]{Lutz Gröll}
\author[First]{Veit Hagenmeyer}

\address[First]{Karlsruhe Institute of Technology, Eggenstein-Leopoldshafen, Baden-Württemberg, Germany (e-mail: anton.ponomarev@kit.edu, lutz.groell@kit.edu, veit.hagenmeyer@kit.edu).}

\begin{abstract}
   We estimate the lock-in domain of the origin of a current control system which is used in common DC/AC inverter designs. The system is a cascade connection of a 4-dimensional linear system (current controller, CC) followed by a two-dimensional nonlinear system (phase-locked loop, PLL). For the PLL, we construct a Lyapunov function via numerical approximation of its level curves. In combination with the quadratic Lyapunov function of the CC, it forms a vector Lyapunov function (VLF) for the overall system. A forward-invariant set of the VLF is found via numerical application of the comparison principle. By LaSalle's invariance principle, convergence to the origin is established.
\end{abstract}


\end{frontmatter}

\section{Introduction}

\subsection{Practical Background}

\emph{Vector current control} is a technique for regulating the AC current in three-phase power electronics. It is based on controlling a vector that represents the three-phase current in a rotating \emph{$dq$ reference frame}. In equilibrium, the reference frame rotates at a constant speed, and the current vector is constant in that frame.

\emph{DC/AC inverters}, also called \emph{grid-tied voltage source converters}, are power electronic devices tasked with supplying DC power to an AC grid. \emph{Grid-following} inverters aim to synchronize with the voltage oscillation that already exists in the grid. To this end, the phase angle of the grid voltage is estimated by the so-called \emph{phase-locked loop} (PLL). The PLL provides a $dq$ frame in which the vector current controller operates.

In the present paper we consider a common inverter design which uses a \emph{synchronous reference frame PLL}, and the current controller itself is a simple \emph{PI controller}.

\subsection{The Problem}

The closed-loop error dynamics of our control system can be described as a \emph{cascade}: the current controller (CC) is a 4-dimensional \emph{linear time-invariant} system that feeds into the PLL which is a two-dimensional \emph{nonlinear} system. The CC subsystem is asymptotically stable, and the PLL resembles a nonlinear pendulum, see \cite[Fig.~1]{ponomarevNonlinearAnalysisSynchronous2024}. If the coordinates of the PLL are $(\dlt\theta, \dlt\omega)$ with $\dlt\theta$ the phase angle error, then the PLL dynamics are $2\pi$-periodic in $\dlt\theta$ with stable equilibria $(2k\pi, 0)$, $k \in \bbZ$, interlaced with unstable ones.

\emph{Asymptotic behavior} of the system is clear: as the CC part exponentially converges to zero, the PLL reduces to an unforced oscillator whose behavior is well documented~-- e.g., see ~\cite{liConditionsExistenceUniqueness2024} and references therein. Under the assumption of small grid impedance, it almost certainly converges to one of the stable equilibria~-- i.e., the set of points \emph{not} converging to one of those has measure zero.

However, the \emph{transient process} is just as important as asymptotic stability: suppose that the system rests at the origin when a sudden disturbance occurs~-- for example, it may be a change in the reference point for the current. This translates to a jump in the CC subsystem. Before the CC returns to zero, it may disturb the PLL enough to cause it to move from $(0,0)$ to another equilibrium like $(2\pi, 0)$. During the transient, $\dlt\theta$ goes through values around $\pi$, i.e., the current is regulated with the phase angle opposite to the correct one. It may lead to power losses, overheating, and may even trigger fault protection.

With the above explanation in mind, it is clear that the \emph{domain of attraction} (DoA) would not be the right concept to describe the behavior of our system since the DoA of the origin includes points with $\dlt\theta$ too far from $0$ which produce undesirable transients.

Instead of the DoA, the notion of a \emph{lock-in domain} has long been in use in the theory of synchronization~-- e.g., \cite{prestonLockinPerformanceAFC1953} analyze it in a frequency control circuit which handles frequency modulation for the purposes of analog color television. A formal definition of the lock-in domain was proposed by \cite{leonovHoldinPullinLockin2015} as the set of initial conditions that do not lead to \emph{cycle slipping} which roughly means that $\dlt\theta$ should not make $2\pi$-turns before settling at an equilibrium.

\emph{The problem} we address in the present paper is to estimate the lock-in domain of the closed-loop system.

\subsection{Existing Results}

There are many results on the stability of grid-tied inverters obtained via \emph{linearization-based analysis} in the frequency domain, e.g., see \cite{vietoSequenceImpedanceModeling2018}. There is also the method of \emph{harmonic linearization} that, unlike the pure linearization, is able to capture coupling between several frequencies, e.g., see \cite{golestanHarmonicLinearizationInvestigation2021}. Such small-signal results are, of course, local, and it is hard to quantify how far from the equilibrium they are applicable.

Among the recent \emph{nonlinear studies}, let us mention \cite{maGeneralizedSwingEquation2022} where the nonlinear transient dynamics of the inverter are investigated using bifurcation theory. An important assumption in that work is \emph{negligible CC dynamics} which is motivated by the claim that the CC is much faster than the PLL. The slow-fast dynamics argument can be formalized by the singular perturbation method of \cite{kokotovicSingularPerturbationMethods1986}~-- this kind of analysis is done by \cite{luSynchronousStabilityAnalysis2024}.

\cite{liConditionsExistenceUniqueness2024} used an energy-like Lyapunov function and LaSalle's invariance principle to estimate the domain of attraction of the zero equilibrium. However, the CC dynamics are neglected in their work as well.

The assumption of negligible CC dynamics is justified in many, but not in all cases. For instance, \cite{midtsundEvaluationCurrentController2010} demonstrate that the interaction between the PLL and CC plays an important role when the inverter is connected to a \emph{weak grid}: the system may become unstable if the grid inductance is too high.

In the present paper, we consider the complete system including not only the fully nonlinear PLL but also the CC dynamics.

\subsection{Our Approach}

This work continues the analysis of the PLL dynamics started in \cite{ponomarevNonlinearAnalysisSynchronous2024}. Therein, we studied the PLL under the influence of a small periodic signal and obtained estimations of a stable periodic solution and its lock-in domain. In the current paper, the PLL is forced by the CC subsystem which produces a signal that is not small but exponentially vanishing. We modify our technique to estimate the lock-in domain of the zero equilibrium. Modification is required to capture the relation between the vanishing disturbance affecting the PLL and the initial conditions of the CC.

Our approach utilizes the method of two-dimensional \emph{comparison systems} which traces back to \cite{elshinMethodComparisonQualitative1954}. A comparison system is a simpler system whose trajectories can only be crossed by the original system in a certain direction. Then the portrait of the original system can in some sense be bounded by the trajectories of the comparison system. \cite{belykhQualitativeInvestigationSecond1975} applied a similar technique to a time-varying system designing an \emph{autonomous} comparison system.

We construct a comparison system for the PLL dynamics treating the CC influence as a bounded time-varying disturbance. This yields nested forward-invariant sets for the PLL. From a collection of such sets we build a Lyapunov function for the PLL. Together with the standard quadratic Lyapunov function of the CC, it forms a two-dimensional vector Lyapunov function (VLF). Applying the comparison principle from the VLF theory of \cite{bellmanVectorLyapunovFunctions1962}, we obtain an invariant set in the phase plane of the VLF which is our final estimation of the lock-in domain.

\subsection{Outline of the Paper}

Section~\ref{se: preliminaries} contains a precise formulation of the system and the goal of the analysis.

In Section~\ref{se: general principles} the Lyapunov functions of the PLL and CC are defined. At this step, we are able to formulate the first estimation of the lock-in domain that we call ``trivial'' (Theorem~\ref{th: trivial estimation}). An improved (larger) estimation is then obtained via the comparison principle from the VLF theory (Theorem~\ref{th: main}).

In Section~\ref{se: realization} we clarify some details concerning the practical realization of the approach. It is explained how the estimations by Theorems~\ref{th: trivial estimation} and~\ref{th: main} can be approximated numerically using standard DAE solvers and optimizers.

Section~\ref{se: numerical} presents a numerical example.

\section{Preliminaries}
\label{se: preliminaries}

\subsection{System Description}

We consider a standard DC/AC inverter current control design based on a synchronous reference frame phase-locked loop and a vector proportional-integral current controller. The output current of the inverter goes through an RL filter. The dynamics of the high-frequency switches that generate the AC voltage are omitted, and the reference point of the current controller is assumed to be constant. The equations of the system can be derived in a standard way, e.g., see~\cite{liConditionsExistenceUniqueness2024} and references therein. We skip the derivation here to save space and only present the result: the \emph{closed-loop error dynamics}
\begin{subequations}
   \label{eq: system}
   \begin{align}
      \label{eq: cc}
      \text{CC:} \quad &
      \hspace{5.75mm} \dot x = Ax, \\
      \label{eq: pll}
      \text{PLL:} \quad &
      \begin{cases}
         \,\dlt\dot\theta = -k_p f(\dlt\theta, \dlt\omega, x) + \dlt\omega, \\
         \dlt\dot\omega = -k_i f(\dlt\theta, \dlt\omega, x)
      \end{cases}
   \end{align}
\end{subequations}
where $k_p > 0$, $k_i > 0$, $\dlt\theta \in \bbR$, $\dlt\omega \in \bbR$, $x \in \bbR^4$, matrix $A$ is Hurwitz,
\begin{equation}
   \label{eq: f}
   f(\dlt\theta, \dlt\omega, x)
   = \frac{g(\dlt\theta, \dlt\omega) - h(\dlt\omega)^T x}
   {\mu - \nu^T x},
\end{equation}
$\mu \in \bbR$ and $\nu \in \bbR^4$ are constant, $g \in \bbR$ and $h \in \bbR^4$ are certain functions with $g(0,0) = 0$ and $\mu \neq 0$. The full expressions of $\mu$, $\nu$, $g$, and $h$ are omitted. They depend on the proportional and integral coefficients $\kappa_p$ and $\kappa_i$ of the current controller which appear in the example (Table~\ref{tb: parameters}).

System~\eqref{eq: system} can be separated into two subsystems which we shall address by name:
\begin{itemize}
   \item subsystem~\eqref{eq: cc} is the \emph{current controller} (CC), written in the synchronous reference frame;
   \item subsystem~\eqref{eq: pll} is the \emph{phase-locked loop} (PLL).
\end{itemize}

Note that system~\eqref{eq: system} has a zero equilibrium. We make the following assumption regarding its stability.

\begin{assumption}
   \label{as: stability}
   The eigenvalues of the Jacobian of~\eqref{eq: pll}
   \begin{equation}
      \begin{bmatrix}
         -\tfrac{k_p}{\mu} f_{\dlt\theta}(0,0,0) &
         \quad 1 - \tfrac{k_p}{\mu} f_{\dlt\omega}(0,0,0) \\[2mm]
         -\tfrac{k_i}{\mu} f_{\dlt\theta}(0,0,0) & \quad 0
      \end{bmatrix}
   \end{equation}
   are non-real with negative real parts (subscripts denote partial derivatives). That is, the zero equilibrium of~\eqref{eq: system} is locally asymptotically stable, and the PLL subsystem~\eqref{eq: pll} is oscillatory near the origin.
\end{assumption}

\begin{rem}
   With Assumption~\ref{as: stability} we limit our attention to the case of \emph{oscillatory} PLL behavior. Non-oscillatory or \emph{overdamped} behavior can be treated similarly.
\end{rem}

\subsection{Goal}

Our goal is to estimate the \emph{lock-in domain} of the zero solution of system~\eqref{eq: system} defined as follows, cf. \cite{ponomarevNonlinearAnalysisSynchronous2024}.

\begin{definition}
   The \emph{lock-in domain} of the origin of system~\eqref{eq: system} is the set of initial conditions from which the solution converges to the origin without crossing the lines $\dlt\theta = \pm\pi$.
\end{definition}

By an \emph{estimation} of the lock-in domain we understand a forward-invariant set contained in the lock-in domain.

\section{General Principles}
\label{se: general principles}

In this section we describe the theoretical approach to the problem. We start with a quadratic Lyapunov function for the CC subsystem, use it to construct a comparison system for the PLL in the sense of \cite{belykhQualitativeInvestigationSecond1975}, and define a Lyapunov function for the PLL via the limit cycles of the comparison system. Then we immediately have a ``trivial'' estimation of the lock-in domain. The estimation is enlarged based on the comparison principle from the vector Lyapunov function theory due to \cite{bellmanVectorLyapunovFunctions1962}.

\subsection{Lyapunov Function for the CC}

As a Lyapunov function of the CC subsystem~\eqref{eq: cc} we take $V\cc(x) = x^T Px$ where $P$ is a positive definite symmetric matrix such that the matrix $A^T P + PA + \gamma P$ is negative definite with some $\gamma > 0$. Then along the solutions of~\eqref{eq: system}
\begin{equation}
   \label{eq: CC Lyapunov derivative}
   \od[V\cc]{t} \leq -\gamma V\cc.
\end{equation}

\subsection{Comparison System for the PLL}

Let us pick a number $V \geq 0$ and consider the PLL subsystem of~\eqref{eq: system} where $x$ is treated as an uncertain time-varying parameter that satisfies $V\cc(x) \leq V$. We apply the general concept of a \emph{left comparison system} \cite[Definition~4]{ponomarevNonlinearAnalysisSynchronous2024} to system~\eqref{eq: system} and arrive at the following definition.

\begin{definition}
   For a given constant $V \geq 0$, system
   \begin{subequations}
      \label{eq: comparison}
      \begin{align}
         \dlt\dot\theta &= -k_p f_*(\dlt\theta, \dlt\omega) + \dlt\omega, \\
         \dlt\dot\omega &= -k_i f_*(\dlt\theta, \dlt\omega)
      \end{align}
      where
      \begin{equation}
         \label{eq: optimization problem}
         f_*(\dlt\theta, \dlt\omega) = \begin{cases}
             \min\limits_{V\cc(x) \leq V} f(\dlt\theta, \dlt\omega, x),
             & \dlt\omega \geq 0, \\
             \max\limits_{V\cc(x) \leq V} f(\dlt\theta, \dlt\omega, x),
             & \dlt\omega < 0
         \end{cases}
      \end{equation}
   \end{subequations}
   is called the \emph{$V$-comparison system} for the PLL~\eqref{eq: pll}.
\end{definition}

\begin{rem}
   The $V$-comparison system is only defined for such $V$ that the set $V\cc(x) \leq V$ does not contain the singularities of $f$, i.e., points with $\nu^T x = \mu$.
\end{rem}

\begin{figure}
   \centering
   \includegraphics{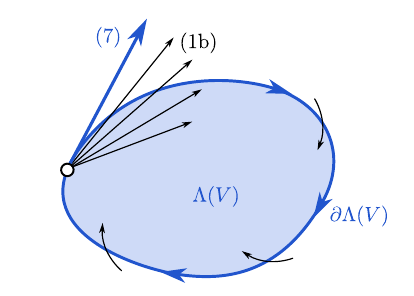}
   \caption{Direction of the $V$-comparison system~\eqref{eq: comparison} is leftmost in the range of possible directions of the PLL subsystem~\eqref{eq: pll} under the constraint $V\cc(x) \leq V$. If a clockwise limit cycle of~\eqref{eq: comparison} bounds a region $\Lambda(V)$ then its boundary $\partial\Lambda(V)$ is crossed by~\eqref{eq: pll} inward.}
   \label{fig: comparison}
\end{figure}

The main property of system~\eqref{eq: comparison} is illustrated by Fig.~\ref{fig: comparison}: under the condition that $V\cc(x) \leq V$, all possible trajectories of~\eqref{eq: pll} cross the trajectories of~\eqref{eq: comparison} rightward, looking in the direction of motion of~\eqref{eq: comparison}. Therefore, the following holds.

\begin{prop}
   \label{pr: limit cycle}
   Suppose that the $V$-comparison system~\eqref{eq: comparison} has a limit cycle orbiting clockwise. Then the inside of the limit cycle is a forward-invariant set of the PLL~\eqref{eq: pll} under the condition $V\cc(x) \leq V$.
\end{prop}

System~\eqref{eq: comparison} is autonomous. For small enough $V$, it is close to~\eqref{eq: pll} with $x = 0$ and thus has a stable limit cycle encompassing the origin clockwise. This follows from the phase portrait of the PLL with $x = 0$ which looks like a nonlinear pendulum~-- e.g., see \cite[Fig.~1]{ponomarevNonlinearAnalysisSynchronous2024}. As $V$ increases, the limit cycle inflates until $V$ reaches a critical value.

Suppose that the limit cycle of~\eqref{eq: comparison} exists for some interval $V \in [0, \bar V]$ (for $V = 0$ it reduces to the zero equilibrium). For each $V \in [0, \bar V]$, let $\Lambda(V)$ be the closed region bounded by the limit cycle corresponding to this value of $V$. The boundary $\partial\Lambda(V)$ is the limit cycle itself. The limit cycles do not intersect each other because that would contradict Proposition~\ref{pr: limit cycle}, i.e.,
\begin{equation}
   \label{eq: matroska}
   \Lambda(V_1) \subset \operatorname{int} \Lambda(V_2) \quad\text{for all}\quad
   0 \leq V_1 < V_2 \leq \bar V
\end{equation}
where $\operatorname{int} \Lambda = \Lambda \setminus \partial\Lambda$.

Let us also note that the limit cycles are smooth because the switching between the minimum and maximum in~\eqref{eq: optimization problem} occurs on the line $\dlt\omega = 0$ where the direction of the vector field of~\eqref{eq: comparison} is independent of the value of $f_*$. That is, although the vector field of~\eqref{eq: comparison} is discontinuous, the limit cycle as a curve in the $(\dlt\theta, \dlt\omega)$-plane is still smooth.

As for how the limit cycle depends on $V$, we make the following assumption.

\begin{assumption}
   The limit cycle $\partial\Lambda(V)$ changes smoothly with $V$ on the interval $[0, \bar V]$.
\end{assumption}

\subsection{Lyapunov Function for the PLL}
\label{se: VPLL}

Let us define the following function $V\pll$ whose sublevel sets are the sets $\Lambda(\cdot)$ described above:
\begin{equation}
   V\pll (\dlt\theta, \dlt\omega) = V \quad\text{for all}\quad
   (\dlt\theta, \dlt\omega) \in \partial\Lambda(V).
\end{equation}
Due to~\eqref{eq: matroska} we have
\begin{equation}
   \label{eq: sublevel sets via VPLL}
   \Lambda(V) = \setb{(\dlt\theta, \dlt\omega): V\pll \leq V}
\end{equation}
and Proposition~\ref{pr: limit cycle} is equivalently formulated as follows.

\begin{lem}
   \label{lem: VPLL decreasing}
   At the points $(\dlt\theta, \dlt\omega, x)$ such that
   \begin{equation}
      V\cc(x) \leq V\pll(\dlt\theta, \dlt\omega)
   \end{equation}
   $V\pll$ is non-increasing along the solutions of system~\eqref{eq: system}.
\end{lem}

\subsection{Trivial Estimation of the Lock-In Domain}

The above construction trivially leads to the following statement.

\begin{thm}
   \label{th: trivial estimation}
   Suppose $\bar V$ is such that the set $\Lambda(\bar V)$ is contained between the lines $\dlt\theta = \pm\pi$. Then the set
   \begin{equation}
      \label{eq: invariant square}
      \setb{
         (\dlt\theta, \dlt\omega, x):
         V\pll(\dlt\theta, \dlt\omega) \leq \bar V, \:
         V\cc(x) \leq \bar V
      }
   \end{equation}
   is an estimation of the lock-in domain of the origin of~\eqref{eq: system}.
\end{thm}

\begin{pf}
   Consider the function
   \begin{equation}
      \label{eq: trivial Lyapunov}
      (\dlt\theta, \dlt\omega, x) \mapsto
      \max\setb{V\pll(\dlt\theta, \dlt\omega), V\cc(x)}.
   \end{equation}
   Suppose that a point $(\dlt\theta, \dlt\omega, x)$ inside~\eqref{eq: invariant square} is such that $V\cc > V\pll$. Then, since $V\cc$ is decreasing by~\eqref{eq: CC Lyapunov derivative}, function~\eqref{eq: trivial Lyapunov} is decreasing as well. On the other hand, if $V\cc \leq V\pll$ then $V\pll$ is non-increasing by Lemma~\ref{lem: VPLL decreasing}, and so is~\eqref{eq: trivial Lyapunov}. Thus, \eqref{eq: trivial Lyapunov} is non-increasing in~\eqref{eq: invariant square}.
   
   Furthermore, the minimal forward-invariant set of~\eqref{eq: system} contained in~\eqref{eq: invariant square} must lie in the plane $x = 0$, and the only such set for the PLL~\eqref{eq: pll} is the origin. Thus, by LaSalle's invariance principle, all solutions inside~\eqref{eq: invariant square} converge to the origin.
\end{pf}

\begin{figure}
   \centering
   \includegraphics{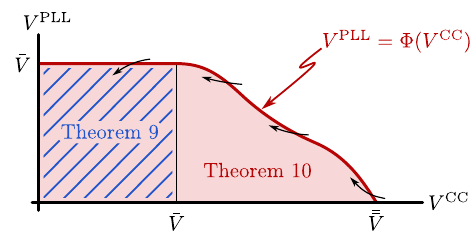}
   \caption{Estimations of the lock-in domain provided by Theorems~\ref{th: trivial estimation} and \ref{th: main} are given in the $(V\cc, V\pll)$-plane. The comparison inequality~\eqref{eq: domain border tangent} prescribes that the boundary should only be crossed inward, as shown by small black arrows.}
   \label{fig: domain}
\end{figure}

\subsection{Improved Estimation}

As the final step, we improve the trivial estimation~\eqref{eq: invariant square} using the \emph{vector Lyapunov function} idea of \cite{bellmanVectorLyapunovFunctions1962}. Consider the vector function
\begin{equation}
   \label{eq: vector Lyapunov}
   (\dlt\theta, \dlt\omega, x) \mapsto \begin{bmatrix*}[l]
      V\cc(x) \\ V\pll(\dlt\theta, \dlt\omega)
   \end{bmatrix*}.
\end{equation}
Whenever $V\pll$ is defined, function~\eqref{eq: vector Lyapunov} maps the state space of~\eqref{eq: system} to the $(V\cc, V\pll)$-plane where the trivial estimation~\eqref{eq: trivial Lyapunov} is the square $[0, \bar V] \times [0, \bar V]$. We would like to extend this square in the direction of larger values of $V\cc$ as shown in Fig.~\ref{fig: domain} so that the estimation would be given by the inequality
\begin{equation}
   \label{eq: domain}
   V\pll \leq \Phi(V\cc)
\end{equation}
with some non-increasing smooth function $\Phi$ such that $\Phi(\cdot) = \bar V$ on $[0, \bar V]$. Function $\Phi$ is defined on some interval $[0, \bar{\bar V}]$ with $\bar{\bar V} > \bar V$ and $\Phi(\bar{\bar V}) = 0$.

To exclude the possibility of system~\eqref{eq: system} crossing the boundary of the set~\eqref{eq: domain} outward, let us request that the \emph{comparison inequality}
\begin{equation}
   \label{eq: domain border tangent}
   \dif V\pll \leq \Phi'(V\cc) \,\dif V\cc
\end{equation}
holds along the solutions of system~\eqref{eq: system} at all points of the boundary, i.e., with $V\pll = \Phi(V\cc)$. Suppose that we know a function $F$ such that
\begin{equation}
   \label{eq: VPLL growth estimation}
   \od[V\pll]{t} \leq F(V\pll, V\cc)
\end{equation}
along~\eqref{eq: system}. Then condition~\eqref{eq: domain border tangent} is satisfied if
\begin{equation}
   F\pb{V\pll, V\cc} \leq \Phi'(V\cc) \od[V\cc]{t}.
\end{equation}
Since $\Phi' \leq 0$ and due to~\eqref{eq: CC Lyapunov derivative}, it is sufficient to demand that the inequality
\begin{equation}
   F\pb{V\pll, V\cc} \leq -\gamma V\cc \Phi'(V\cc)
\end{equation}
holds on the boundary $V\pll = \Phi(V\cc)$. This leads to the following statement.

\begin{thm}
   \label{th: main}
   Suppose $\bar V$ is such that the set $\Lambda(\bar V)$ is contained between the lines $\dlt\theta = \pm\pi$. Let $\Phi(V)$ be a function given by
   \begin{equation}
      \Phi(V) = \bar V \quad\text{for all}\quad V\in[0, \bar V]
   \end{equation}
   and as a solution to the differential equation
   \begin{equation}
      \label{eq: ode Phi}
      \Phi'(V) = -\frac{F\pb{\Phi(V), V}}{\gamma V}, \quad
      \Phi(\bar V) = \bar V
   \end{equation}
   on the interval $[\bar V, \bar{\bar V}]$ with $\bar{\bar V}$ such that $\Phi(\bar{\bar V}) = 0$ where $F$ is a function that satisfies~\eqref{eq: VPLL growth estimation} along the solutions of~\eqref{eq: system}. Then the set
   \begin{equation}
      \label{eq: domain in coordinates}
      \setb{(\dlt\theta, \dlt\omega, x) :
      V\pll(\dlt\theta, \dlt\omega) \leq \Phi\pb{V\cc(x)}}
   \end{equation}
   is an estimation of the lock-in domain of the origin of~\eqref{eq: system}.
\end{thm}

\begin{pf}
   By the comparison argument preceding the theorem, set~\eqref{eq: domain in coordinates} is forward-invariant. Since $V\cc$ is strictly decreasing, the system eventually enters the set~\eqref{eq: invariant square} where Theorem~\ref{th: trivial estimation} asserts convergence to the origin.
\end{pf}

\section{Realization of the Method}
\label{se: realization}

In order to carry out the construction of Theorems~\ref{th: trivial estimation} and~\ref{th: main} numerically, one needs to complete the following:
\begin{itemize}
   \item simulate the $V$-comparison system~\eqref{eq: comparison} and compute its limit cycles for different values of $V$, thus approximating the function $V\pll$;
   \item approximate function $F$ defined by~\eqref{eq: VPLL growth estimation} and solve the ODE~\eqref{eq: ode Phi} for $\Phi$.
\end{itemize}
For these tasks, we suggest the following tools.

\subsection{Solving the Comparison System}

To find an optimality condition for the problem~\eqref{eq: optimization problem}, we compute the gradient of $f$ with respect to $x$:
\begin{equation}
   \label{eq: nabla f}
   \nabla_x f = \frac{b(\dlt\theta, \dlt\omega) + Hx}
   {(\mu - \nu^T x)^2}
\end{equation}
where $b = \mu h(\dlt\omega) - g(\dlt\theta, \dlt\omega) \nu \in \bbR^4$ and $H = h(\dlt\omega) \nu^T - \nu h(\dlt\omega)^T = \const \in \bbR^{4\times 4}$. We observe that the gradient is never zero. Therefore, the extremum in~\eqref{eq: optimization problem} is reached on the boundary, and optimality conditions can be written as
\begin{subequations}
   \begin{align}
      x^T P x &= V, \\
      Px &= \lambda \nabla_x f(\dlt\theta, \dlt\omega, x)
   \end{align}
\end{subequations}
where Lagrange multiplier $\lambda$ is positive in case of maximization and negative for minimization. We combine the optimality conditions together with the dynamics of~\eqref{eq: comparison} in a single system of differential-algebraic equations (DAE)
\begin{subequations}
   \label{eq: comparison dae}
   \begin{align}
      \dlt\dot\theta &= -k_p f_*(\dlt\theta, \dlt\omega) + \dlt\omega, \\
      \dlt\dot\omega &= -k_i f_*(\dlt\theta, \dlt\omega), \\
      0 &= x\mmin^T P x\mmin - V, \\
      0 &= Px\mmin - \lambda\mmin \pb{
         b(\dlt\theta, \dlt\omega) + Hx\mmin
      }, \\
      0 &= x\mmax^T P x\mmax - V, \\
      0 &= Px\mmax - \lambda\mmax \pb{
         b(\dlt\theta, \dlt\omega) + Hx\mmax
      }
   \end{align}
   where
   \begin{equation}
      f_*(\dlt\theta, \dlt\omega) = \begin{cases}
         f(\dlt\theta, \dlt\omega, x\mmin), & \dlt\omega \geq 0, \\
         f(\dlt\theta, \dlt\omega, x\mmax), & \dlt\omega < 0.
      \end{cases}
   \end{equation}
\end{subequations}
This is a 12-dimensional system where the state contains $\dlt\theta$ and $\dlt\omega$ as well as two 4-vectors $x\mmin$ and $x\mmax$ which are the minimizer and maximizer of problem~\eqref{eq: optimization problem} and two scalars $\lambda\mmin < 0$ and $\lambda\mmax > 0$ which are the corresponding Lagrange multipliers. System~\eqref{eq: comparison dae} has index 1 and can be tackled, e.g., by the DAE solvers included with Matlab.

\begin{rem}
   The DAE approach relies on the fact that the minimizers $x\mmin$ and $x\mmax$ in our problem~\eqref{eq: optimization problem} depend at least piecewise continuously on $\dlt\theta$ and $\dlt\omega$.
\end{rem}

\subsection{Approximating the Function $F$}

To satisfy~\eqref{eq: VPLL growth estimation}, we are looking for the function
\begin{multline}
   \label{eq: F as optimization}
   F(\tilde V\pll, \tilde V\cc) = \max_{\substack{
      V\pll(\dlt\theta, \dlt\omega) = \tilde V\pll \\
      V\cc(x) \leq \tilde V\cc
   }} \\
   \nabla V\pll(\dlt\theta, \dlt\omega)
   \cdot \begin{bmatrix*}[l]
      -k_p f(\dlt\theta, \dlt\omega, x) + \dlt\omega \\
      -k_i f(\dlt\theta, \dlt\omega, x)
   \end{bmatrix*}
\end{multline}
where ``$\cdot$'' is the dot product. Assuming we know the gradient $\nabla V\pll$, optimization~\eqref{eq: F as optimization} can be done approximately by iterating through the level curves of $V\pll$.

The problem is now to find $\nabla V\pll$. Observe that $\nabla V\pll$ is related to the linear part of the variation of the limit cycle $(\dlt\theta, \dlt\omega)$ of system~\eqref{eq: comparison dae} with respect to the parameter $V$. To find the linear part, we differentiate~\eqref{eq: comparison dae} with respect to $V$:
\begin{subequations}
   \label{eq: gradient dae}
   \begin{align}
      \dlt\dot\theta' &= -k_p f_*' + \dlt\omega', \\
      \dlt\dot\omega' &= -k_i f_*', \\
      0 &= 2x\mmin^T P x\mmin' - 1, \\
      0 &= Px\mmin' - \lambda\mmin' \pb{
         b(\dlt\theta, \dlt\omega) + Hx\mmin
      } \notag\\
      &\phantom{={}} - \lambda\mmin \p{
         b_{\dlt\theta} \dlt\theta'
         + b_{\dlt\omega} \dlt\omega'
         + H x\mmin'
      }, \\
      0 &= 2x\mmax^T P x\mmax' - 1, \\
      0 &= Px\mmax' - \lambda\mmax' \pb{
         b(\dlt\theta, \dlt\omega) + Hx\mmax
      } \notag\\
      &\phantom{={}} - \lambda\mmax \p{
         b_{\dlt\theta} \dlt\theta'
         + b_{\dlt\omega} \dlt\omega'
         + H x\mmax'
      }
   \end{align}
   where
   \begin{equation}
      f_*' = \begin{cases}
         f_{\dlt\theta} \dlt\theta' + f_{\dlt\omega} \dlt\omega'
         + (\nabla_x f)^T x\mmin', & \dlt\omega \geq 0, \\
         f_{\dlt\theta} \dlt\theta' + f_{\dlt\omega} \dlt\omega'
         + (\nabla_x f)^T x\mmax', & \dlt\omega < 0,
      \end{cases}
   \end{equation}
\end{subequations}
prime denotes the derivative with respect to $V$, and subscripts denote partial derivatives. System~\eqref{eq: gradient dae} is another 12-dimensional system of DAEs whose state is $(\dlt\theta', \dlt\omega', \lambda\mmin', x\mmin', \lambda\mmax', x\mmax')$.

Suppose that the combined system~\eqref{eq: comparison dae}, \eqref{eq: gradient dae} with a certain value of $V$ has a limit cycle
\begin{equation}
   \pb{
      \dlt\theta(t), \dlt\omega(t), \dots, \dlt\theta'(t), \dlt\omega'(t), \dots
   }.
\end{equation}
Then the fact $V\pll(\dlt\theta, \dlt\omega) = V$ implies two equations
\begin{equation}
   \label{eq: equations for gradient}
   \nabla V\pll \cdot \begin{bmatrix}
      \dlt\dot\theta \\ \dlt\dot\omega
   \end{bmatrix} = 0 \quad\text{and}\quad
   \nabla V\pll \cdot \begin{bmatrix}
      \dlt\theta' \\ \dlt\omega'
   \end{bmatrix} = 1
\end{equation}
from which $\nabla V\pll$ is determined along the limit cycle.

\subsection{A Sketch of the Implementation}

Systems~\eqref{eq: comparison dae} and~\eqref{eq: gradient dae} together form a 24-dimensional index-1 system of DAEs that we solve in Matlab while increasing $V$ in small steps, at each step waiting until the system converges close enough to a limit cycle. The process is stopped once the system escapes the range $-\pi \leq \dlt\theta \leq \pi$. The final value of $V$ is used as $\bar V$. The collection of thus obtained limit cycles is an approximation of the level sets of $V\pll$. Knowing $\dlt\theta'$ and $\dlt\omega'$ along the limit cycles, we approximate $\nabla V\pll$ by solving two linear equations~\eqref{eq: equations for gradient}. Function $F$ is found via approximate optimization in~\eqref{eq: F as optimization} along the limit cycles $(\dlt\theta, \dlt\omega)$ and interpolation between them. Equation~\eqref{eq: ode Phi} is then integrated to attain function $\Phi$ that bounds the estimation~\eqref{eq: domain in coordinates}.

\begin{table}[b]
   \centering
   \caption{Parameters of the numerical example}
   \begin{tabular}{lll}
      \multicolumn{3}{c}{Inverter} \\
      \hline
      & version I & version II \\
      $k_p$ & $3 \times 10^{-4}$ & \,$3 \times 10^{-3}$ \\
      $k_i$ & $1 \times 10^{-4}$ & \,$1 \times 10^{-2}$ \\
      $\kappa_p$ & \multicolumn{2}{l}{\hspace{7mm} $1 \times 10^{-2}$} \\
      $\kappa_i$ & \multicolumn{2}{l}{\hspace{7mm} $1$} \\
      $L\filter$ & \multicolumn{2}{l}{\hspace{7mm} $1 \times 10^{-3}$} \\
      $R\filter$ & \multicolumn{2}{l}{\hspace{7mm} $4 \times 10^{-4}$} \\
      $i\dqref$ & \multicolumn{2}{l}{\hspace{7mm} $[10 \:\:\: 0]^T$}
   \end{tabular}
   \hspace{4mm}
   \begin{tabular}{ll}
      \multicolumn{2}{c}{Grid} \\
      \hline
      $\omega\grid$ & $100\pi$ \\
      $L\grid$ & $2 \times 10^{-3}$ \\
      $R\grid$ & $6 \times 10^{-4}$ \\
      $\norm{v\grid}$ & 325 \\
      \\
      \\
      \\
      \\
   \end{tabular}
   \label{tb: parameters}
\end{table}

\section{Numerical Results}
\label{se: numerical}

Parameters used for the numerical example are shown in Table~\ref{tb: parameters}. There are two versions of the PLL coefficients $k_p$ and $k_i$ that we are going to compare.

\begin{figure}[t]
   \centering
   \includegraphics[trim = 6.5cm 11.1cm 6.5cm 11.7cm, clip]{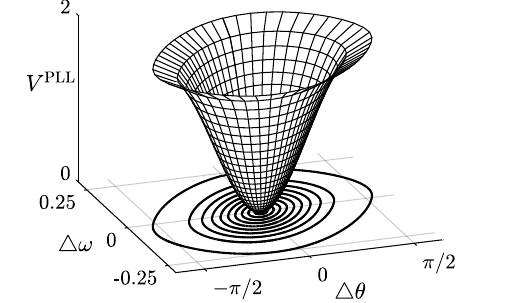}
   \caption{Surface plot and level curves of the Lyapunov function $V\pll$ of the PLL subsystem~\eqref{eq: pll} with version~I of parameters $k_p$ and $k_i$ in the numerical example.
   }
   \label{fig: example levels}
\end{figure}

For version~I, Fig.~\ref{fig: example levels} shows the Lyapunov function $V\pll$ of the PLL subsystem~\eqref{eq: pll} constructed according to the idea of Section~\ref{se: VPLL} and method of Section~\ref{se: realization}. The largest sublevel set of $V\pll$ is interpretable as an estimation of the lock-in domain of the PLL assuming there is no error in the current ($V\cc = 0$)~-- this set extends to $\dlt\theta \approx \pm\pi/2$ which is far enough to be meaningful. For version~II of $k_p$ and $k_i$, function $V\pll$ looks similar and is not shown.

For both versions~I and II, Fig.~\ref{fig: example domain} presents the main result: an estimation of the lock-in domain in the $(V\cc, V\pll)$-plane. To explain the difference between the two versions, we note that version~II is a 10 times faster PLL than version~I. Within the hatched region (the \emph{trivial estimation} given by Theorem~\ref{th: trivial estimation}), $V\pll$ strictly decreases for both versions, and the speed of the decrease is mostly irrelevant, hence the trivial regions are almost the same for both versions. However, the extended region (Theorem~\ref{th: main}) for version~II is much smaller~-- this is because here $V\pll$ can \emph{increase} temporarily until $V\cc$ drops enough that the system reaches the trivial region. The faster version~II is prone to ``overreacting'' and leaving the lock-in domain~-- a trait that can be predicted directly from the PLL equations~\eqref{eq: pll}: increasing $k_p$ and $k_i$ makes the dynamics faster but at the same time more sensitive to the disturbance $x$.

\section{Conclusion}

\begin{figure}[t!]
   \centering
   \includegraphics[trim = 6.5cm 11.8cm 6.5cm 12.55cm, clip]{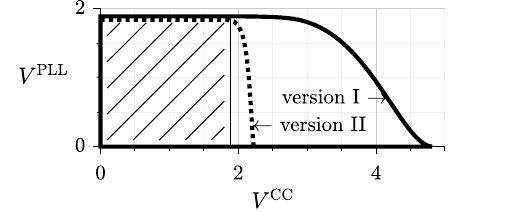}
   \caption{Estimation of the lock-in domain of system~\eqref{eq: system} by Theorem~\ref{th: main} in the numerical example with two versions of the PLL parameters. The hatched part is the trivial estimation by Theorem~\ref{th: trivial estimation} (almost the same for both versions).}
   \label{fig: example domain}
\end{figure}

We introduce a method to estimate the lock-in domain of the origin of the current control system~\eqref{eq: system}. The estimation is given as an inequality involving two Lyapunov functions $V\cc$ and $V\pll$ designed for the subsystems~\eqref{eq: cc} and~\eqref{eq: pll}. It has the shape shown in Fig.~\ref{fig: example domain}, i.e., a square $[0, \bar V] \times [0, \bar V]$ with a bulge on its side. In our implementation, we tried finding the largest possible $\bar V$. It would be possible to take smaller values of $\bar V$ and thus obtain smaller invariant sets attracted to the origin. A collection of such sets would define a Lyapunov function of the system near the origin.

The technique can also be directly applied to other systems structured as a linear subsystem feeding into a nonlinear two-dimensional subsystem.


\bibliography{ifacconf}

\begin{thebibliography}{13}
\providecommand{\natexlab}[1]{#1}
\providecommand{\url}[1]{\texttt{#1}}
\providecommand{\urlprefix}{URL }
\expandafter\ifx\csname urlstyle\endcsname\relax
  \providecommand{\doi}[1]{doi:\discretionary{}{}{}#1}\else
  \providecommand{\doi}{doi:\discretionary{}{}{}\begingroup
  \urlstyle{rm}\Url}\fi

\bibitem[{Bellman(1962)}]{bellmanVectorLyapunovFunctions1962}
Bellman, R. (1962).
\newblock Vector {{Lyapunov}} functions.
\newblock \emph{Journal of the Society for Industrial and Applied Mathematics
  Series A Control}, 1(1), 32--34.
\newblock \doi{10.1137/0301003}.

\bibitem[{Belykh(1975)}]{belykhQualitativeInvestigationSecond1975}
Belykh, V.N. (1975).
\newblock The qualitative investigation of a second order nonautonomous
  nonlinear equation.
\newblock \emph{Differentsial'nye Uravneniya}, 11(10), 1738--1753.
\newblock In Russian.

\bibitem[{Elshin(1954)}]{elshinMethodComparisonQualitative1954}
Elshin, M.I. (1954).
\newblock The method of comparison in the qualitative theory of an incomplete
  differential equation of second order.
\newblock \emph{Mat. Sbornik N.S.}, 34/76, 323--330.
\newblock In Russian.

\bibitem[{Golestan et~al.(2021)Golestan, Guerrero, Vasquez, Abusorrah, and
  Al-Turki}]{golestanHarmonicLinearizationInvestigation2021}
Golestan, S., Guerrero, J.M., Vasquez, J.C., Abusorrah, A.M., and Al-Turki, Y.
  (2021).
\newblock Harmonic linearization and investigation of three-phase
  parallel-structured signal decomposition algorithms in grid-connected
  applications.
\newblock \emph{IEEE Transactions on Power Electronics}, 36(4), 4198--4213.
\newblock \doi{10.1109/TPEL.2020.3021723}.

\bibitem[{Kokotovi{\'c} et~al.(1986)Kokotovi{\'c}, Khalil, and
  O'Reilly}]{kokotovicSingularPerturbationMethods1986}
Kokotovi{\'c}, P.V., Khalil, H.K., and O'Reilly, J. (1986).
\newblock \emph{Singular perturbation methods in control: analysis and design}.
\newblock Academic Press, London.

\bibitem[{Leonov et~al.(2015)Leonov, Kuznetsov, Yuldashev, and
  Yuldashev}]{leonovHoldinPullinLockin2015}
Leonov, G.A., Kuznetsov, N.V., Yuldashev, M.V., and Yuldashev, R.V. (2015).
\newblock Hold-in, pull-in, and lock-in ranges of {{PLL}} circuits: Rigorous
  mathematical definitions and limitations of classical theory.
\newblock \emph{IEEE Transactions on Circuits and Systems I: Regular Papers},
  62(10), 2454--2464.
\newblock \doi{10.1109/TCSI.2015.2476295}.

\bibitem[{Li et~al.(2024)Li, Lu, Tang, and
  Du}]{liConditionsExistenceUniqueness2024}
Li, Y., Lu, Y., Tang, Y., and Du, Z. (2024).
\newblock Conditions of existence and uniqueness of limit cycle for
  grid-connected {{VSC}} with {{PLL}}.
\newblock \emph{IEEE Transactions on Power Systems}, 39(1), 706--719.
\newblock \doi{10.1109/TPWRS.2023.3238000}.

\bibitem[{Lu et~al.(2024)Lu, Liu, Shao, Lin, Tian, Yang, Li, and
  Du}]{luSynchronousStabilityAnalysis2024}
Lu, Y., Liu, M., Shao, C., Lin, J., Tian, G., Yang, Q., Li, Y., and Du, Z.
  (2024).
\newblock Synchronous stability analysis model of grid-following converter
  based on singular perturbation method.
\newblock In \emph{2024 10th {{International Conference}} on {{Electrical
  Engineering}}, {{Control}} and {{Robotics}} ({{EECR}})}, 343--348.
\newblock \doi{10.1109/EECR60807.2024.10607206}.

\bibitem[{Ma et~al.(2022)Ma, Li, Kurths, Cheng, and
  Zhan}]{maGeneralizedSwingEquation2022}
Ma, R., Li, J., Kurths, J., Cheng, S., and Zhan, M. (2022).
\newblock Generalized swing equation and transient synchronous stability with
  {{PLL-based VSC}}.
\newblock \emph{IEEE Transactions on Energy Conversion}, 37(2), 1428--1441.
\newblock \doi{10.1109/TEC.2021.3137806}.

\bibitem[{Midtsund et~al.(2010)Midtsund, Suul, and
  Undeland}]{midtsundEvaluationCurrentController2010}
Midtsund, T., Suul, J.A., and Undeland, T. (2010).
\newblock Evaluation of current controller performance and stability for
  voltage source converters connected to a weak grid.
\newblock In \emph{The 2nd {{International Symposium}} on {{Power Electronics}}
  for {{Distributed Generation Systems}}}, 382--388.
\newblock \doi{10.1109/PEDG.2010.5545794}.

\bibitem[{Ponomarev et~al.(2024)Ponomarev, Hagenmeyer, and
  Gröll}]{ponomarevNonlinearAnalysisSynchronous2024}
Ponomarev, A., Hagenmeyer, V., and Gröll, L. (2024).
\newblock Nonlinear analysis of the synchronous reference frame phase-locked
  loop under unbalanced grid voltage.
\newblock \emph{Nonlinear Dynamics}, 112, 9225--9243.
\newblock \doi{10.1007/s11071-024-09532-9}.

\bibitem[{Preston and Tellier(1953)}]{prestonLockinPerformanceAFC1953}
Preston, G.W. and Tellier, J.C. (1953).
\newblock The lock-in performance of an {{AFC}} circuit.
\newblock \emph{Proceedings of the IRE}, 41(2), 249--251.
\newblock \doi{10.1109/JRPROC.1953.274214}.

\bibitem[{Vieto and Sun(2018)}]{vietoSequenceImpedanceModeling2018}
Vieto, I. and Sun, J. (2018).
\newblock Sequence impedance modeling and analysis of {{Type-III}} wind
  turbines.
\newblock \emph{IEEE Transactions on Energy Conversion}, 33(2), 537--545.
\newblock \doi{10.1109/TEC.2017.2763585}.

\end{thebibliography}

\end{document}